\documentclass[a4paper,11pt]{article}
\usepackage{amsthm,amsmath}
\usepackage{dsfont}
\usepackage{mathptmx}
\usepackage{algorithm,algpseudocode}
\usepackage{tikz}
\usepackage{pgfplots}
\usepackage{cite}
\usepackage{color}
\usepackage{latexsym}
\usepackage{lscape}
\usepackage[latin1]{inputenc}
\usepackage{graphicx}
\usepackage{booktabs}
\usepackage[T1]{fontenc}
\usepackage{authblk}
\usepackage{url}
\usepackage{caption}
\usepackage{subcaption}
\usepackage{todonotes}

\newcommand{\be}{\begin{equation}}
	\newcommand{\ee}{\end{equation}}

\newcommand{\bc}{\begin{center}}
	\newcommand{\ec}{\end{center}}

\title{Comparative numerical results on HANSO (Hybrid Algorithm for Nonsmooth Optimization) }
\author[1]{Ali Hakan TOR}
\affil[1]{Abdullah G\"{u}l University, Computer Science Faculty, Departemnt of Applied Mathematics}

\begin{document}
\maketitle

\begin{abstract}
The purpose of this study is to give an opinion on which software should be used by researchers working in applied fields when using free software found in the literature. Of course, this study is not enough alone as there is no comparison with other free software, but this study provides comprehensive information about the algorithms used by the Hybrid Algorithm for Non-smooth Optimization software, known as HANSO. In addition to this information, the robustness and the accuracy of this software are discussed in this study. 
\end{abstract}


\smallskip

\section{Introduction}
In this study, some numerical results have been obtained by using the software HANSO (Hybrid Algorithm for Non-Smooth Optimization). HANSO is a hybrid version of the Broyden, Fletcher, Goldfarb, and Shanno Algorithm (BFGS) and the Gradient Sampling Algorithm (GSA). If you run HANSO by default, it gives you the minimum value of an optimization problem by using the BFGS and GSA together. On the other hand, if you want, it also offers the opportunity to run BFGS and GSA algorithms separately. In this study, we will obtain results from the HANSO and BFGS algorithms written by \cite{HANSO}, and we will observe how HANSO works well as a hybrid method. The reason we compare HANSO with BFGS instead of GSA is that HANSO uses BFGS for the most part of its calculations by seeing GSA only as an auxiliary algorithm. You can find the BFGS ang GSA algorithms coded by Overton in the HANSO package on the website \cite{HANSO}. In addition to all these, I would like to emphasize that the BFGS implementation we use for comparison is the same as the BFGS implementation used by HANSO. Therefore, we can make a fair observation to understand how much HANSO hybridizes BFGS with GSA. 

The organization of this study proceeds as follows. In the following section, the software HANSO is examined by explaining how it works and which algorithms have been used. Then, a review for these algorithms is given with details. In the next section, the characteristic of the test problems which are used for numerical experience are given briefly. All numerical results are reported and discussed in the final section. 

\section{HANSO (Hybrid Algorithm for Non-Smooth Optimization)}

Version 2.2 of HANSO developed by Michael Overton is used in this study. It has General Public License (GNU) as published by the Free Software Foundation, so anybody can redistribute it and/or modify it under the terms of this license. 

HANSO is intended to seek a minimum value of non-smooth, non-convex functions, but also applicable to functions that are smooth, convex or both. It is based on the BFGS algorithm and GSA. You can find some details about BFGS and GSA in the following subsections.

\subsection{BFGS Algorithm}

BFGS algorithm suggested independently by Broyden, Fletcher, Goldfarb, and Shanno, in 1970 uses the Quasi-Newton algorithm which is a generalization of the secant method. The main difference between BFGS and Quasi-Newton algorithms is that it uses and maintains different properties of the matrix when updating formulas. In BFGS, the Hessian matrix is not calculated. Instead of this calculation, BFGS uses a inverse Hessian matrix approximation using information from gradient evaluation. BFGS is normally used for optimizing smooth, not necessarily convex, functions, for which the convergence rate is generically superlinear. However, BFGS has acceptable performance even for non-smooth optimization problems, typically with a linear convergence rate as long as a weak Wolfe line search is used. This version of BFGS will work well both for smooth and non-smooth functions and has a stopping criterion that applies for both cases \cite{lewis12}.

There are several options for the stopping criterion of BFGS algorithm in HANSO. First of all, it is possible to adjust the tolerance of a decent direction. If its norm is less than the given tolerance, the code is terminated. In this study, the default tolerance $10^{-6}$ is used. Another stopping criterion is that the distance of the gradient vector calculated in each step from the current iteration point is greater than the given tolerance value. The default tolerance value $10^{-4}$ is used in this study again. Other stopping criteria are related to change of function values, the magnitude of the current iteration point and CPU time, but in this study numerical results have been obtained without any restrictions on them.

\subsection{Gradient Sampling Algorithm (GSA)}
Gradient sampling idea was used in \cite{ermo82,shor85} for the first time. Later, the gradient sampling method was used to approximate the Clarke subdifferential for locally Lipschitz functions in \cite{burk021} and it was
improved for non-smooth non-convex problems in \cite{burk05}, which is used in HANSO Software. Later, other versions of gradient sampling methods for some special optimization problems was developed such as \cite{burk06,burk04,kiwiel07}.

GSA is intended for non-convex  and locally Lipschitz functions that are differentiable almost everywhere, in other words, they are not differentiable on a set of the measure zero, so the subgradient at a randomly selected point is uniquely determined as the gradient at that point. Therefore, in GSA, gradients are computed on a set of randomly generated nearby points at current iteration. Consequently, by using gradient sampling, a local information of the function is obtained and the quadratic subproblem is formed. The $\epsilon$-steepest descent direction is constructed by solving this quadratic subproblem, where $\epsilon$ is the sample radii.

The stopping criterion of GSA in HANSO is on descent directions. If the norm of the descent direction at current iteration is less then given tolerance, the algorithm is terminated. HANSO's default values $10^{-6}$ is used as a tolerance in this study.

\section{Test Problem}
The efficiency of HANSO was tested on the well-known nonsmooth optimization academic test problems taken in \cite{luksan00}. The reasons why all test problems in \cite{luksan00} were not included are different. The first reason is the repetition of some problems, namely, CB2 and Rosen-Suzuki, in both Chapters 2 and 3 of \cite{luksan00}. The second reason is the unboundedness of some problems, that is there is no global solution to these problems. The names of these problems are Bard, Gamma, Colville 1 and HS78. After that, several problems, namely, PBC3, Kowalik-Osborne, EXP, PBC1, EVD61, and Filter, have more than one local solution. Another reason is that the input data are not fully available for the problem TR48, so the problem TR48 is not placed. Lastly, the problem Transformer has complex coefficients, so it is not used. Briefly, we use 36 test problems from both Chapters 2 and 3 in \cite{luksan00}, while there are 50 test problems. While some of them are nonconvex, all these test problems are nonsmooth (for detail see Table \ref{tableofproblems}). 
\begin{table}[h!]
\begin{center}
	\begin{tabular}{ c c c c}
		\textbf{Problem} & \textbf{Number of Variable} & \textbf{Optimal Value} & \textbf{Convexity} \\ 
		\hline
		CB2 & 2 & 1,9522245& Convex \\  
		WF & 2 & 0 & Nonconvex\\
		SPIRAL&2&0&Nonconvex\\
		Rosenbrock& 2 & 0 & Nonconvex\\
		Crescent& 2& 0& Nonconvex\\
		CB3& 2 & 2& Convex \\
		DEM & 2 & -3 & Convex\\
		QL & 2 & 7,2 & Convex\\
		LQ & 2 & -1,4142136& Convex\\
		Mifflin 1 & 2 & -1 & Convex\\
		Mifflin 2 & 2 & -1 & Nonconvex\\
		Wolfe & 2 & -8 & Convex \\
		EVD52 & 3 & 3,5997193& Nonconvex\\ 
		Rosen$\_$Suzuki & 4 & -44 & Convex\\
		Polak6 & 4 & -44 & Convex\\
		Davidon 2 & 4 & 115,70644& Convex\\
		OET5 & 4 & $0,26359735 \times 10^{-2}$ & Nonconvex\\
		OET6 & 4 & $0,20160753 \times 10^{-2}$ & Nonconvex\\
		Shor & 5 & 22,600162 & Convex\\
		El-Attar & 6&0,5598131 & Nonconvex\\
		Wong 1 & 7& 680,63006& Convex\\
		Wong 2 & 10 & 24,306209 & Convex\\
		Polak 2 & 10 & 54,598150& Convex \\
		Maxquad & 10 & -0,8414083 & Convex\\
		Gill & 10 & 9,7857721 & Nonconvex \\
		Polak 3 & 11 & 3,70348& Convex\\
		Osborne 2 & 11 & $0,48027401 \times 10^{-1}$ & Nonconvex\\
		Steiner 2 & 12 & 16,703838& Nonconvex \\
		Shell Dual & 15 & 32,348679 & Nonconvex\\
		Wong 3 & 20 & 93,90525& Convex \\
		Watson & 20 & $0,14743027 \times 10^{-7}$ & Convex\\
		Maxq & 20 & 0 & Convex\\
		Maxl & 20 & 0& Convex\\
		Gofflin & 50 & 0 & Convex\\
		MXHILB & 50 & 0 & Convex\\
		L1HILB & 50 & 0 & Convex\\
	\end{tabular}
\end{center}
\caption{List of Academic Test Problems}
\label{tableofproblems}
\end{table}

HANSO allows us to use the specified starting point or randomly generated starting point. However, in this study, 20 randomly generated starting points were used. Since the HANSO code was not suitable for running 20 different points, this code was modified to use these randomly generated 20 points by reading our starting point files. These numerical results are presented in the next section.

\section{Numerical Results}

Since it is not possible to give all of the 20 results obtained for each of the 36 test problems mentioned in the previous section, we first give the table below, which presents HANSO and BFGS solves how many problems related to the starting points successfully.

\begin{table}[h!]
	\begin{center}
		\begin{tabular}{ c c c c}
			\textbf{Problem} & \textbf{HANSO} & \textbf{BFGS} & \textbf{Convexity} \\ 
			\hline
			CB2 & 20 & 20& Convex \\  
			WF & 2 & 2 & Nonconvex\\
			SPIRAL&18&18&Nonconvex\\
			Rosenbrock& 20 & 20 & Nonconvex\\
			Crescent& 0& 0& Nonconvex\\
			CB3& 14 & 14& Convex \\
			DEM & 1 & 1 & Convex\\
			QL & 20 & 20 & Convex\\
			LQ & 20 & 20& Convex\\
			Mifflin 1 & 0 & 0 & Convex\\
			Mifflin 2 & 0 & 0 & Nonconvex\\
			Wolfe & 20 & 20 & Convex \\
			EVD52 & 20 & 20& Nonconvex\\ 
			Rosen$\_$Suzuki & 0 & 0 & Convex\\
			Polak6 & 0 & 0 & Convex\\
			Davidon 2 & 20 & 20& Convex\\
			OET5 & 0 & 0 & Nonconvex\\
			OET6 & 0 & 0 & Nonconvex\\
			Shor & 20 & 20 & Convex\\
			El-Attar & 8 & 8 & Nonconvex\\
			Wong 1 & 0 & 0 & Convex\\
			Wong 2 & 20 & 20 & Convex\\
			Polak 2 & 0 & 0 & Convex \\
			Maxquad & 20 & 20 & Convex\\
			Gill & 3 & 3 & Nonconvex \\
			Polak 3 & 5 & 5 & Convex\\
			Osborne 2 & 0 & 0 & Nonconvex\\
			Steiner 2 & 0 & 0 & Nonconvex \\
			Shell Dual & 0 & 0 & Nonconvex\\
			Wong 3 & 10 & 10 & Convex \\
			Watson & 0 & 0 & Convex\\
			Maxq & 20 & 20 & Convex\\
			Maxl & 20 & 20 & Convex\\
			Gofflin & 20 & 20 & Convex\\
			MXHILB & 20 & 20 & Convex\\
			L1HILB & 20 & 20 & Convex\\
			\hline
			\textbf{TOTAL} & 361 & 361 &  
		\end{tabular}
	\end{center}
	\caption{Number of problems solved successfully}
	\label{sucsolv}
\end{table}

We can observe that there is no difference in the number of solving problems between HANSO and BFGS from Table \ref{sucsolv}. In addition, when there were 20 starting points for 36 test problems, in other words, these algorithms worked 720 times, they were successful in only $\% 50$ of these works. This means that the number of iterations, function and gradient evaluations for these softwares should be compared for a more effective comparison.

We will make these comparisons for two different situations. First, without looking at whether the related software solves the problem or not, we will compare the average number of iterations, functions and gradients evaluation for each problem. Then, we will do the same comparison only for the problems it solved. I mean, for example, WF was solved 2 times successfully, so the average values of these two successful solutions will be compared. 

Table \ref{allaverage} presents the average number of iterations, functions and gradients evaluations among the results obtained running with 20 starting points for each problem. In Table \ref{allaverage}, $n_{iter}$, $n_{fev}$ and $n_{gev}$ denote the numbers of iterations, functions and gradients, respectively.

\begin{table}[h!]
	\begin{center}
		\begin{tabular}{ c| c c c |c c c}
			\textbf{Problem} & \multicolumn{3}{|c|}{\textbf{HANSO}} & \multicolumn{3}{c}{\textbf{BFGS}} \\ 
		\hline
			& $n_{iter}$ & $n_{fev}$ & $n_{gev}$ & $n_{iter}$ & $n_{fev}$ & $n_{gev}$\\
			\hline
			CB2 & 27,9 & 61,3& 32,9 & 35,2 & 74,3& 61,3\\  
			WF & 6 & 40,1 & 7,4 & 6,5 & 41,2 & 9,4\\
			SPIRAL&353,1&608,8&353,2&357,6&618,9&365,8\\
			Rosenbrock& 40,2 & 79 & 45,2& 40,2 & 79 & 45,2\\
			Crescent& 1,6& 2& 1,6 & 1,6& 2& 1,6\\
			CB3& 22,4 & 66,7& 26,2& 57,4 & 173,2& 166,4 \\
			DEM & 3,2 & 34 & 4,4 & 3,2 & 34 & 4,4\\
			QL & 23,1 & 444,9 & 654 & 30,2 & 85 & 53,3\\
			LQ & 15,3 & 34,6& 16,3& 39,3 & 112,4& 111,5\\
			Mifflin 1 & 12,6 & 48 & 13,6 & 43,4 & 198,8 & 135,8\\
			Mifflin 2 & 19,3 & 87,7 & 20,3& 19,6 & 89,4 & 21,4\\
			Wolfe & 24,6 & 73,2 & 39,5 & 29,6 & 83,1 & 59,5 \\
			EVD52 & 46,3 & 137,4& 62,3 & 58,5 & 175,8 & 133,1\\ 
			Rosen$\_$Suzuki & 16,1 & 804,7 & 953,1& 16,1 & 71,8 & 17,1\\
			Polak6 & 16,7 & 114 & 21,9& 16,7 & 114 & 21,9\\
			Davidon 2 & 156,8 & 350& 359& 205,4 & 491,4& 747,1\\
			OET5 & 905 & 2374,5 & 2994,5& 987,9 & 2591,6 & 3459\\
			OET6 & 311,9 & 833,8 & 850,5& 337,1 & 889,1 & 1023\\
			Shor & 44,6 & 109 & 102,1& 58,8 & 142,7 & 242,3\\
			El-Attar & 513,1 & 1350,4 & 1688,6& 513,1 & 1350,4 & 1688,6\\
			Wong 1 & 21,5 & 96,4 & 22,1& 21,5 & 96,4 & 22,1\\
			Wong 2 & 275,1 & 690,2 & 1666,7& 316,7 & 805,1 & 2499,7\\
			Polak 2 & 12,9 & 64 & 19,5& 12,9 & 64 & 19,5 \\
			Maxquad & 215,4 & 543,3 & 1303,3& 215,4 & 543,3 & 1303,3\\
			Gill & 1000 & 1697,1 & 2708,7& 1000 & 1697,1 & 2708,7 \\
			Polak 3 & 612,5 & 1644,2 & 1796,1& 612,5 & 1650,3 & 1825,2\\
			Osborne 2 & 511,1 & 1468,6 & 1382,5& 623,3 & 1794,9 & 1568,3\\
			Steiner 2 & 25,9 & 77,9 & 26,6& 25,9 & 77,9 & 26,6 \\
			Shell Dual & 3,9 & 47,5 & 4,4& 3,9 & 47,5 & 4,4\\
			Wong 3 & 477,9 & 1152 & 942,9& 477,9 & 1152 & 942,9 \\
			Watson & 2,4 & 36,6 & 2,9& 2,4 & 36,6 & 2,9\\
			Maxq & 208,9 & 413,5 & 520,5& 208,9 & 413,5 & 520,5\\
			Maxl & 1000 & 3896,5 & 26428& 1000 & 3896,5 & 26428\\
			Gofflin & 645,4 & 2713,9 & 20197& 645,4 & 2713,9 & 20197\\
			MXHILB & 1000 & 5152,1 & 47514,1& 1000 & 5152,1 & 47514,1\\
			L1HILB & 821,9 & 1787,7 & 1629,9& 821,9 & 1787,7 & 1629,9\\
 
		\end{tabular}
	\end{center}
	\caption{The Average Numbers}
	\label{allaverage}
\end{table}

A very striking result in Table \ref{allaverage} is immediately noticeable. For 20 problems, there is no difference between HANSO and  BFGS. This means that although HANSO is called a hybrid method, it only uses the GSA algorithm for 16 problems. When we look at 16 problems that give different results for HANSO and BFGS, we see that HANSO achieved its goal and decreased $n_{iter}$, $n_{fev}$ and $n_{gev}$ except for the problem Rosen-Suzuki. In addition to these, in problems where the average number of iterations is 1000, the software terminated due to the restriction on the number of iterations.

Finally, we will make a comparison between the average numbers for all starting points and the average numbers for the starting points at which the problem solved. In this comparison, 8 problems, namely WF, SPIRAL, CB3, DEM, El-Attar, Gill, Polak 3 and Wong 3  are used. Other problems are not used because they either have solved for all starting points or not solved for any starting point. There is no difference to compare in both case. In Table \ref{8average}, for these 8 problem, the average numbers $n_{iter}$, $n_{fev}$ and $n_{gev}$ for all starting points and for the successfully solved starting points are given in the first row and the second row of the related problem, respectively.

\begin{table}[h!]
	\begin{center}
		\begin{tabular}{ c| c c c |c c c}
			\textbf{Problem} & \multicolumn{3}{|c|}{\textbf{HANSO}} & \multicolumn{3}{c}{\textbf{BFGS}} \\ 
			\hline
			& $n_{iter}$ & $n_{fev}$ & $n_{gev}$ & $n_{iter}$ & $n_{fev}$ & $n_{gev}$\\
			\hline
			WF & 6 & 40,1 & 7,4 & 6,5 & 41,2 & 9,4\\
			 & 17 & 48,5& 23,5 & 22 & 59& 43\\
			 \hline  
			SPIRAL&353,1&608,8&353,2&357,6&618,9&365,8\\
			& 289 & 496,6 & 289,2& 294,2 & 508,4 & 304,2\\
			\hline  
			CB3& 22,4 & 66,7& 26,2& 57,4 & 173,2& 166,4 \\
			 & 27,6 & 68,8 & 32,6 & 77,6 & 220,9 & 232,9\\
			 \hline  
			DEM & 3,2 & 34 & 4,4 & 3,2 & 34 & 4,4\\
			 & 26 & 53& 30& 26 & 53& 30\\
			 \hline  
			El-Attar & 513,1 & 1350,4 & 1688,6& 513,1 & 1350,4 & 1688,6\\
			 & 321,9 & 784,8 & 1367,5& 321,9 & 784,8 & 1367,5\\
			 \hline  
			Gill & 1000 & 1697,1 & 2708,7& 1000 & 1697,1 & 2708,7 \\
			 & 1000 & 1815,3 & 1146,3& 1000 & 1815,3 & 1146,3\\
			 \hline  
			Polak 3 & 612,5 & 1644,2 & 1796,1& 612,5 & 1650,3 & 1825,2\\
			 & 807,8 & 2151,2 & 2609,2& 814,8 & 2175,8 & 2725,6 \\
			 \hline  
			Wong 3 & 477,9 & 1152 & 942,9& 477,9 & 1152 & 942,9 \\
			 & 597,2 & 1449,1 & 1525,2& 597,2 & 1449,1 & 1525,2\\
		
		\end{tabular}
	\end{center}
	\caption{The Average Numbers}
	\label{8average}
\end{table}

From Table 4, we see that the average values increase for the starting points at which the problems solved in these 8 problems except Spirial, El-Atar and Gill. This means that both HANSO and BFGS need more iteration for successfully solved problems than for unsolved problems. In Spiral and El-Attar, the situation is a little different, the softwares have made more calculations at the starting points where they could not solve the problem. In Gill, there is no change in the number of iterations, since the maximum number of iterations is reached. It can be easily observed that the number of function and gradient evaluation have not shoved a similar trend. 

\section{Conclusion}

When we consider all these numerical results, we can say that HANSO, which is a hybrid method, does not give very different results from BFGS in terms of accuracy and correctness, but it reduces the number of iterations and calculations to some extent. On the other hand, when we look at the results obtained with different starting points, it can be said that for many problems, both HANSO and BFGS are sensitive to the starting point, that is, neither of them is robust. Nevertheless, this free software can be used for academic studies. Of course, when using this software, it would be more appropriate to use more than one starting point instead of one starting point for academic study. Finally, there is no need to use HANSO and BFGS separately to get results. The results obtained with only HANSO are sufficient.

\end{document}